\documentclass[reqno,12pt,centertags]{article}
\usepackage{amsmath,amsthm,amscd,amssymb,latexsym,upref}
\date{\today}

\input epsf
\usepackage{epsfig}
\usepackage[T2A,OT1]{fontenc}
\usepackage[ot2enc]{inputenc}
\usepackage[russian,english]{babel}

\usepackage{epic,eepicemu}

\newcommand{\bbD}{{\mathbb{D}}}

\newcommand{\bbZ}{{\mathbb{Z}}}

\newcommand{\cH}{{\mathcal{H}}}
\newcommand{\bbT}{{\mathbb{T}}}

\newcommand{\cF}{{\mathcal{F}}}
\newcommand{\cL}{{\mathcal{L}}}

\newcommand{\cE}{{\mathcal{E}}}

\newcommand{\cK}{{\mathcal{K}}}
\newcommand{\cM}{{\mathcal{M}}}
\newcommand{\cN}{{\mathcal{N}}}

\newcommand{\fe}{{\mathfrak{e}}}
\newcommand{\ff}{{\mathfrak{f}}}

\newcommand{\fA}{{\mathfrak{A}}}

\renewcommand{\Re}{\text{\rm Re}}

\newcommand{\tr}{\text{\rm tr}}

\newcommand{\Sz}{{\mathbf{Sz}}}
\newcommand{\HSz}{{\mathbf{HS}}}
\newcommand{\GI}{{\mathbf{GI}}}

\allowdisplaybreaks \numberwithin{equation}{section}
\newtheorem{theorem}{Theorem}[section]

\newtheorem{lemma}[theorem]{Lemma}
\newtheorem{proposition}[theorem]{Proposition}

\theoremstyle{definition}
\newtheorem{definition}[theorem]{Definition}

\newtheorem{remark}[theorem]{Remark}
\newtheorem{problem}[theorem]{Problem}

\title
{Faddeev-Marchenko scattering for CMV matrices and the Strong Szeg\"o Theorem}
\author{L. Golinskii, A. Kheifets, F. Peherstorfer\thanks{
Partially supported by the Austrian Founds FWF, project number:
P20413--N18}, and P. Yuditskii$^*$}

\date{\today}
\begin{document}
\maketitle
\begin{abstract}
Simon  proved the existence of the wave operators for the CMV matrices with Szeg\"o class Verblunsky coefficients, and therefore the existence of the scattering function.
Generally, there is no hope to restore a CMV matrix when we start  from the scattering function, in particular, because it does not contain any information about the (possible) singular measure. 
 Our main point of interest is the solution of the inverse scattering problem (the heart of the Faddeev--Marchenko theory),
 that is, to give necessary and sufficient conditions on a certain class of CMV matrices  such that the restriction of this correspondence (from a matrix to the scattering function)  is one to one.
 In this paper we show that the main questions on inverse scattering can be solved with the help of
three important classical results: Adamyan-Arov-Krein (AAK) Theory, Helson-Szeg\"o Theorem and Strong Szeg\"o Limit Theorem. Each of these theorem states the equivalence  of certain conditions. Actually, to each theorem we add one more equivalent condition related to the CMV inverse scattering problem.
\end{abstract}
\section{Introduction}

To a given collection of numbers $\{a_n\}_{n\ge 0}$ in  the unit open disk $\bbD$ and $a_{-1}$ in the unit circle $\bbT$
 we associate the CMV matrix $\fA=\fA_1\fA_0$, where
\begin{equation}\label{matrixform0}
    \fA_0=\begin{bmatrix}
    A_0 & & \\ &
     A_2 & \\
    & & \ddots
    \end{bmatrix},\fA_1= \begin{bmatrix}
-\bar a_{-1}& & & \\ &
A_1& & \\
   & & A_3 & \\
  &  & & \ddots
    \end{bmatrix}
\end{equation}
and the $A_k$'s are the $2\times 2$ unitary matrices
$$
A_k=\begin{bmatrix}
 a_k&\rho_k\\ \rho_k&-\bar a_k
    \end{bmatrix},\quad \rho_k=\sqrt{1-|a_k|^2}.
     $$

Note that $\fA$ is a unitary operator in $l^2(\bbZ_+)$. The initial vector $e_0$ of the standard basis is cyclic for $\fA$, indeed by the definition
\begin{equation}\label{gi4}
\begin{split}
    \fA^{-1}\{\rho_{2n-1}e_{2n-1}-e_{2n}\bar a_{2n-1}\}=&e_{2n}\bar
a_{2n}+e_{2n+1}\rho_{2n}\\
\fA\{e_{2n}\rho_{2n}-e_{2n+1} a_{2n}\}=&e_{2n+1}
a_{2n+1}+e_{2n+2}\rho_{2n+1}.
\end{split}
\end{equation}
That is, acting in turn by $\fA^{-1}$ and $\fA$ on $e_0$ we can get in the linear combination any vector of the standard basis.

Let $\sigma$ be the spectral measure of $\fA$, i.e.,
\begin{equation}\label{gi5}
   R(z)= \int_{\bbT}\frac{t+z}{t-z}d\sigma=\langle\frac{\fA+z}{\fA-z}e_0,e_0\rangle.
\end{equation}

The matrix $\fA$ is of the Szeg\"o class, $\fA\in \Sz$, if $\sum |a_k|^2<\infty$.
It holds if and only if
 the measure $\sigma$ is of the form
$$
d\sigma=w(t) dm(t)+d\sigma_{s},\quad \log w\in L^1,
$$
where $dm(t)$ is the Lebesgue measure and $\sigma_{s}$ is the singular component. Define the outer function
$$
D(z)=e^{\frac 1 2\int_{\bbT}\frac{t+z}{t-z}\log w(t) dm(t)}, \quad D_*(z)=\overline{D(1/\bar z)}.
$$

Simon \cite[Sect. 10]{Sv2} proved the existence of the wave operators for the CMV matrices $\fA$ with $l^2$ (Szeg\"o class) Verblunsky coefficients, and therefore the existence of the scattering function, which is of the form
\begin{equation}\label{desssim}
    s(t)=-a_{-1}\frac{D(t)}{D_*(t)}.
\end{equation}
Generally, there is no hope to restore $\fA$ when we start  from $s(t)$, in particular, because it does not contain any information about the (possible) singular measure. Even if we assume that
 $\fA\in \Sz_{a.c.}$, i.e., $\sigma_s=0$, the correspondence $\fA\mapsto s$ is not one to one on its image.
(An easy example: $s(t)=t^2$ corresponds simultaneously to $D(t)=(1-t)^2$ and $D(t)=(1+t)^2$).
 Our main point of interest is the solution of the inverse scattering problem (the heart of the Faddeev--Marchenko theory),
 that is, to give necessary and sufficient conditions on a certain class of CMV matrices and correspondingly on the associated scattering functions such that the restriction of the map $\fA\mapsto s$ is one to one. Naturally, we would like to have an explicit algorithm to get $\fA$. One of the  key elements of the Faddeev--Marchenko construction is the so called Gelfand-Levitan-Marchenko (GLM) transformation operators, which transform an orthogonal standard basis into an intrinsic orthogonal basis related to a perturbed CMV matrix.

Usually in the Faddeev-Marchenko theory the scattering function appears in the following way.

\begin{theorem}\label{th1}
Let $\fA\in \Sz$. Then there exists a unique generalized eigenvector
$\Psi(t)=\{\Psi_n(t)\}_{n=0}^\infty$ and a so called scattering function $s(t)$, $|s(t)|=1$, such that
\begin{equation}\label{gi2}
    \begin{bmatrix}\Psi_0(t)&\Psi_1(t)&\dots\end{bmatrix}\fA=
   t \begin{bmatrix}\Psi_0(t)&\Psi_1(t)&\dots\end{bmatrix},\ t\in \bbT,
\end{equation}
and the following asymptotics are satisfied
\begin{equation}\label{gi3}
    \Psi_{2n}(t)=t^n+o(1),\quad \Psi_{2n+1}(t)=\overline{s(t)} t^{-n-1}+o(1), \ n\to \infty,
\end{equation}
in the $L^2$--norm.
\end{theorem}

Also,
$s(t)$ is called the nonlinear Fourier transform (NLFT) of the Verblunsky coefficients $\{a_k\}_{k=-1}^\infty$ \cite{TTh}.

In fact, Theorem \ref{th1} is a restatement of the classical Szeg\"o theorem on the asymptotic behavior of orthogonal polynomials on the unit circle (OPUC).

 Note, $\fA$ is unitary equivalent to the multiplication operator by the independent variable
   in $L^2_{d\sigma}$ with respect to the following orthonormal basis, see \eqref{gi4},
\begin{equation}\label{gi6}
\begin{split}
    t^{-1}\{\rho_{2n-1}P_{2n-1}(t)-P_{2n}(t)\bar a_{2n-1}\}=&P_{2n}(t)\bar
a_{2n}+P_{2n+1}(t)\rho_{2n}\\
t\{P_{2n}(t)\rho_{2n}-P_{2n+1}(t) a_{2n}\}=&P_{2n+1}(t)
a_{2n+1}+P_{2n+2}(t)\rho_{2n+1},
\end{split}
\end{equation}
where $P_n$ are Laurent polynomials
$$
P_0(t)=\pi^{(0)}_0,\ P_1(t)=\pi^{(1)}_1 t^{-1}+\pi^{(1)}_0,\ P_2(t)=\pi^{(2)}_2 t+\pi^{(2)}_1 t^{-1}+\pi^{(2)}_0,\dots
$$

Now we can see that the relations
\begin{equation}\label{gi8}
    \lim_{n\to\infty}t^{-n}{D_*(t)}P_{2n}(t)=1,\quad
    \lim_{n\to\infty}t^{n+1}{D(t)}P_{2n+1}(t)=-\bar a_{-1}
\end{equation}
are indeed consequences of Szeg\"o's Theorem. Moreover, $\Psi_n(t)={D_*(t)} P_n(t)$,
 which  proves simultaneously the representation \eqref{desssim} of $s$ in \eqref{gi3}.


In this paper we show that the above posed questions on inverse scattering can be solved with the help of
three important classical results: Adamyan-Arov-Krein (AAK) Theory, Helson-Szeg\"o Theorem and Strong Szeg\"o Limit Theorem. Each of these theorem states the equivalence  of certain conditions. Actually, to each theorem we add one more equivalent condition related to the CMV inverse scattering problem. In this way we describe three subclasses of $\Sz_{a.c.}$:
\begin{itemize}
\item
$\Sz^{reg}\subset \Sz_{a.c.}$ on which $\fA\mapsto s$ is one to one,
\item
$\fA\in \HSz\subset \Sz^{reg}$ if the corresponding  GLM operator $\cM$ is bounded,
\item
$\fA\in \GI\subset \HSz$ if $\det(\cM^*\cM)<\infty$.
\end{itemize}


The AAK theory describes the solutions of the so called Nehari problem, which are of the form
\begin{equation*}\label{aak100}
  -\frac{\psi}{\bar\psi}\frac{\cE+\bar\phi}{1+\phi\cE},
\end{equation*}
where $\phi$ depends on the data of the problem;
 $\psi$, $\psi(0)>0$, is an outer function, such  that $|\psi|^2+|\phi|^2=1$;  and
 $\cE$ is an arbitrary function of the Schur class ($\cE\in H^\infty, \|\cE\|\le 1)$. The function
 $\phi$ has the following properties
 \begin{equation}\label{propphi100}
    \phi\in H^\infty,\ \|\phi\|\le 1,\
    \phi(0)=0,\ \log(1-|\phi|^2)\in L^1.
 \end{equation}
 But not any function of the form \eqref{propphi100} generates the set of \textit{all} solutions of the Nehari problem. If $\phi$ generates all solutions it is called regular. The AAK theory gives several necessary and sufficient conditions for regularity. Our contribution to them is as follows:
 Let $R$ be given by \eqref{gi5} and let
 $$
 \phi(z)=a_{-1}\frac{R(0)-R(z)}{R(0)+R(z)}.
 $$
 Then $\fA\mapsto s$ is one to one if and only if $\phi$ is regular.

 A unimodular function $s$ belongs to the Helson-Szeg\"o class, $s\in HS$, if it possesses the representation
\begin{equation}\label{hsconds}
    s=e^{i(\tilde u+v)},\  u,v\in L^\infty, \ \|v\|<\pi/2,
\end{equation}
where $\tilde u$ is the harmonic function conjugated to $u$.
The classical   Helson-Szeg\"o and Hunt-Muckenhoupt-Wheeden Theorems give several necessary and sufficient conditions for a function $s$ to be of Helson--Szeg\"o class. Our contribution is the following:
$\fA\mapsto s$ is one to one and the GLM operator is bounded if and only if $s\in HS$.

The B.Golinskii--Ibragimov Theorem is discussed in details in Sect. 6 of Simon's book \cite{Sv1}. Our Theorem complements this theory from the point of view of  scattering. Recall   a function $f(t)$ belongs to the Sobolev space $B_2^{1/2}$ if its Fourier coefficients $\{c_k\}_{k=-\infty}^\infty$
satisfy the condition
$$
\sum_k|k| |c_k|^2<\infty.
$$
In a sense Theorem \ref{th5.3} says that the Nonlinear Fourier Transform (NLFT) of the Verblunsky coefficients (that is the scattering function $s(t)$) belongs to $B_2^{1/2}$ if and only if its Linear Fourier Transform (LFT) belongs to $B_2^{1/2}$. We believe that this Theorem is an essential improvement of the general Faddeev-Marchenko theory (see  Remark \ref{rk5.5}).

Finally we demonstrate  that the so-called Widom Formula has a natural proof in the frame of the scattering theory.

\section{From spectral representation to scattering representation}
First of all let us point out that \eqref{gi5} gives a  one to one correspondence between
the CMV matrix $\fA$ on one hand and the measure $\sigma$ and normalizing constant $a_{-1}$ on the other hand.
In the orthogonalization procedure for the system
$$
1, t^{-1}, t, t^{-2}, t^{2}, t^{-3},\dots
$$
with respect to $d\sigma$ we choose
$$
\pi^{(0)}_0>0, \quad \frac{\pi^{(1)}_1}{|\pi^{(1)}_1|}=-\bar a_{-1}.
$$
In this case
\begin{equation}\label{pirho}
 \pi^{(2n)}_{2n} =\frac 1{\rho_0\rho_1\dots\rho_{2n-1}},\quad
  \pi^{(2n+1)}_{2n+1} =-\bar a_{-1}\frac 1{\rho_0\rho_1\dots\rho_{2n}}.
\end{equation}

In the Szeg\"o case we have the decomposition $d\sigma=|D|^2dm+\sigma_s$ and we define the scattering function $s=-a_{-1}D/D_*$, thus we have the map
\begin{equation}\label{cor}
    \fA\leftrightarrow\{d\sigma, a_{-1}\}\mapsto s.
\end{equation}

To clarify the uniqueness property of this map we define certain Hilbert  spaces associated with the symbol $s$. But first we give an another description of the space $L^2_{d\sigma}$.
For $f\in L^2_{d\sigma}$ we set
\begin{equation}\label{toscat}
    \cF(z)=\cF(z,f)=(F_1(z), F_2(z)),
\end{equation}
where
\begin{equation}\label{toscat1}
\begin{split}
F_1(z)=F_1(z;f)=&\frac 1{D(z)}\int_{\bbT}\frac{t}{t-z}f(t)d\sigma(t),
\quad |z|<1,\\
F_2(z)=F_2(z;f)=&\frac{a_{-1}}{D_*(z)}\int_{\bbT}\frac{t}{t-z}f(t)d\sigma(t),
\quad |z|>1.
\end{split}
\end{equation}
The first component is analytic inside  the unit disk, the second in its exterior. Actually, they are functions of bounded characteristic with an outer denominator (of the Smirnov class).

Note that for $z\in \bbD$
\begin{equation}\label{tomsigma}
    (DF_1)(z)-\bar a_{-1}(D_*F_2)(1/\bar z)=\int_{\bbT}\frac{1-|z|^2}{|t-z|^2}f(t) d\sigma(t).
\end{equation}

\begin{definition} Let  $\cF=(F_1(z),F_2(z))$, where $F_1(z)$ is  analytic in $\bbD$ and $F_2(z)$, $F_2(\infty)=0$, is analytic for $|z|>1$. We say that $\cF\in M(\sigma,a_{-1})$ if
\begin{equation}\label{defmsig}
 \|\cF\|^2=  \sup_{r<1}\int_{\bbT}|(DF_1)(rt)-\bar a_{-1}(D_*F_2)(\bar t/r)|^2\frac{dm(t)}{\Re R(rt)}<\infty,
\end{equation}
where $R(z)$ is defined in \eqref{gi5}.
\end{definition}

\begin{theorem} $\cF\in M(\sigma, a_{-1})$ if and only if it is of the form
\eqref{toscat1}. Moreover
\begin{equation}\label{normnorm}
    \|\cF\|_{M(\sigma, a_{-1})}=\|f\|_{L^2_{d\sigma}}.
\end{equation}
\end{theorem}
\begin{proof} The proof is based on a concept of the so called Hellinger integral,
see e.g. \cite{shm} or \cite{KhYu, Kats}. We recall the corresponding construction here.

Let
\begin{equation}\label{hell1}
\begin{bmatrix}    da &db\\
* &dc \end{bmatrix}
\end{equation}
be $2\times 2$ nonnegative  matrix measure on $\bbT$.
In fact, it means that $db$ is absolutely continuous with respect to $da$,
$db=fda$, moreover
$$
dc-|f|^2da \ge 0.
$$
Let us point out that $\inf c(\bbT)$ over all possible nonnegative matrix-functions
 of the form \eqref{hell1} with the fixed $da$ and $db$ corresponds precisely to the case
  \begin{equation}\label{hell2}
    dc=|f|^2da.
\end{equation}

Now, let
\begin{equation*}\label{hell3}
    \begin{bmatrix}    u(z) & v(z)\\
* & w(z) \end{bmatrix}=\int_{\bbT}\frac{1-|z|^2}{|t-z|^2}\begin{bmatrix}    da(t) &db(t)\\
* &dc(t) \end{bmatrix}.
\end{equation*}
Since this matrix is nonnegative we have
\begin{equation*}\label{hell4}
    w_r(t)\ge \frac{|v_r(t)|^2}{u_r(t)},\quad t\in\bbT,\ r\in(0,1),
\end{equation*}
where $u_r(z)=u(rz)$, etc.
Therefore
\begin{equation}\label{hell5}
    c(\bbT)=w(0)\ge \int_\bbT\frac{|v_r(t)|^2}{u_r(t)} dm(t).
\end{equation}

Thus we get
\begin{equation}\label{hell6}
    c(\bbT) \ge \sup_r I(r),\ \text{where}\ I(r)=\int\frac{|v_r(t)|^2}{u_r(t)} dm(t),
\end{equation}
for all $dc$ which forms a nonnegative matrix function \eqref{hell1} with the given $da$ and $db$.
Note that this already proves the easy part of the theorem. Indeed, set
\begin{equation}\label{hell7}
    u(z)=\Re R(z), \quad v(z)=(DF_1)(z)-\bar a_{-1}(D_*F_2)(1/\bar z).
\end{equation}
Due to \eqref{tomsigma} we have boundedness of $\sup$ in \eqref{defmsig}.

Now, for harmonic functions $u(z)$ and $v(z)$ let us prove that $I(r)$ increases with $r$.
For $k\in(0,1)$, the matrix measure
\begin{equation*}\label{hell8}
    \begin{bmatrix}    u_k(t) & v_k(t)\\
* & \frac{|v_k(t)|^2}{u_k(t)} \end{bmatrix}dm(t)
\end{equation*}
is nonnegative  on $\bbT$.
Its harmonic extension in the disk is of the form
\begin{equation}\label{hell9}
    \begin{bmatrix}    u_k(z) & v_k(z)\\
* & w(z,k) \end{bmatrix}.
\end{equation}
Therefore, due to \eqref{hell5} we get
\begin{equation*}\label{hell10}
    I(k)=w(0,k)\ge \int \frac{|v_k(rt)|^2}{u_k(rt)} dm(t)=\int \frac{|v(rkt)|^2}{u(rkt)} dm(t)=I(rk).
\end{equation*}

Thus, if $\sup I(r)<\infty$ then the limit $\lim_{r\to1}I(r)$ exists.
Consider the sequence of harmonic matrix functions \eqref{hell9}. Due to the compactness principle there exists a limit (on a subsequence)
\begin{equation*}\label{hell11}
   \lim_{k_n\to 1}\begin{bmatrix}    u_{k_n}(z) & v_{k_n}(z)\\
* & w(z,k_n) \end{bmatrix}= \begin{bmatrix}    u(z) & v(z)\\
* & w(z,1) \end{bmatrix}=\int_{\bbT}\frac{1-|z|^2}{|t-z|^2}\begin{bmatrix}    da(t) &db(t)\\
* &dc(t,1) \end{bmatrix}.
\end{equation*}
 Therefore, for this particular matrix measure  we have $db(t)=f(t)da(t)$, and
moreover,
$$
c(\bbT,1)=w(0,1)=\lim I(k_n)=\sup I(r).
$$
 Since generally we have  inequality \eqref{hell6}, this value of $c(\bbT,1)$ corresponds to the infimum, and, therefore $dc(t,1)$ is of the form \eqref{hell2}, that is,
$$
\int|f|^2da=\sup I(r).
$$
Thus for $u$ and $v$ of the form \eqref{hell7} we get \eqref{normnorm}.
\end{proof}

\begin{lemma} $M(\sigma, a_{-1})$ is a space with the reproducing kernels $\cK_z=\cK_z(\sigma, a_{-1})$,
$|z|<1, |z|>1$, in particular,
$$
\langle \cF,  \cK_0\rangle= F_1(0),\quad \langle \cF,  \cK_\infty\rangle= (tF_2)(\infty),
$$
where
\begin{equation}\label{rk8}
\begin{split}
    \cK_0=&\left(\frac{1}{D D(0)}\frac{R(0)+R}{2},\frac{a_{-1}}{D_* D(0)}\frac{R(0)-\bar R}{2}\right),\\
    \cK_\infty=&\left(\frac{\bar a_{-1}}{D D(0)}\frac{R(0)-R}{2t},\frac{1}{D_* D(0)}\frac{R(0)+\bar R}{2t}\right).
    \end{split}
\end{equation}
\end{lemma}

\begin{proof} Use definitions \eqref{toscat1} and \eqref{normnorm}.
\end{proof}

Now we note the following identity  for the boundary values
$$
{|D(t)|^2}{f(t)}=(DF_1)(t)-\bar a_{-1}(D_*F_2)(t),\quad t\in\bbT.
$$
We introduce the scalar product
\begin{equation}\label{ts1}
    \|\cF\|_{s}^2:=\int|s(t)F_1(t)+F_2(t)|^2 dm(t)=\int|f|^2|D|^2dm\le\int|f|^2d\sigma.
\end{equation}
Note that equality holds if and only if $f|_{supp(\sigma_s)}=0$.

\begin{definition}
Let $M_s=M_s(D)$ be the Hilbert space of the functions \eqref{toscat}, $f|_{supp(\sigma_s)}=0$, with the scalar product
\eqref{ts1}.
\end{definition}
Let us point out that the scalar product depends on $s$ but the collection of functions depends actually on $D$, that is why, it's better to write $M_s(D)$.

\begin{definition}
We define $\check M_s=
\text{clos}\{\cF=(F_1,F_2): F_1\in H^2,\ F_2\in H^2_-\}$ with the norm,
$$
\|\cF\|_s^2=\int|s(t)F_1(t)+F_2(t)|^2 dm(t)=\langle\begin{bmatrix} I&\cH_s^*\\\cH_s& I\end{bmatrix}
\begin{bmatrix} F_1\\ F_2\end{bmatrix},\begin{bmatrix} F_1\\ F_2\end{bmatrix}\rangle,
$$
where $\cH=\cH_s$ is the Hankel operator with the symbol $s$
$$
\cH_s: H^2\to H^2_-, \quad \cH_s F_1= P_-(sF_1),\ F_1\in H^2.
$$
\end{definition}
The space $\check M_s$ depends only on $s$, more precisely on $\cH_s$, that is, of the negative Fourier coefficients of $s$, so we write $\check M_s(\cH)$.
\begin{lemma} Let $D\in H^2$, $D(0)>0$,  be an outer function and $s=-a_{-1}D/\bar D$, $a_{-1}\in\bbT$. Then
 $\check M_s(\cH_s)\subset M_s(D)$.
\end{lemma}

\begin{proof} For polynomials $P_1,P_2$ set
$
f=\frac 1{\bar D}P_1+\frac 1 D \overline{tP_2}\in L^2_{wdm}.
$
Then $F_1(f)=P_1, F_2(f)=-a_{-1}\overline{tP_2}$. This set is complete in $\check M_s$.
\end{proof}


\section{AAK Theory and uniqueness in Inverse Scattering }

The AAK Theory deals with the Nehari problem \cite{AAK1, AAK2, AAK3}, see also \cite{Garnett}.

\begin{problem} Given a Hankel operator $\cH$, $\|\cH\|\le 1$, describe the collection of symbols
 \begin{equation}\label{nehclass}
    \cN(\cH)=\{f\in L^\infty: \cH=\cH_f,\ \|f\|_\infty\le 1\}.
\end{equation}
 \end{problem}
 The Nehari Theorem stays solvability of the problem, i.e., $\cN(\cH)\not=\emptyset$.
The following lemma is related to the question of uniqueness of a solution.
\begin{lemma}[Adamyan-Arov-Krein]
The point evaluation functional
$$\cF\to F_1(0)$$
is bounded in
 $\check M_s(\cH)$ if and only if
 \begin{equation}\label{aak1}
    \lim_{r\uparrow 1}\langle(I-r^2\cH^*\cH)^{-1}1,1\rangle<\infty.
 \end{equation}
 Moreover,
 \begin{equation}\label{defchk}
    \check \cK_0=\check \cK_0^s=\lim_{r\uparrow 1}\begin{bmatrix} I&r\cH^*\\r\cH& I\end{bmatrix}^{-1}\begin{bmatrix} 1\\ 0\end{bmatrix},\ \
\check \cK_\infty=\check \cK_\infty^s=\lim_{r\uparrow 1}\begin{bmatrix} I&r\cH^*\\r\cH& I\end{bmatrix}^{-1}\begin{bmatrix} 0\\ \bar t\end{bmatrix}.
 \end{equation}
\end{lemma}

\begin{theorem}[Adamyan-Arov-Krein] A solution of the Nehari problem is not unique if and only if
\eqref{aak1} holds. In this case
the set $\cN(\cH)$ is of the form
\begin{equation}\label{dcrp-nehclass}
    \cN(\cH)=\{f=f_{\cE}=-\frac{\psi_\cH}{\bar\psi_\cH}\frac{\cE+\bar\phi_\cH}{1+\phi_\cH\cE}: \cE\in
 H^\infty,\ \|\cE\|_\infty\le 1\},
\end{equation}
where $\phi_\cH$ is a Schur class function given by
\begin{equation}\label{defphi}
    \phi_\cH(z)=\frac{\overline{(\check{\cK}_0)_2(1/\bar z)}}{(\check{\cK}_0)_1(z)}=z\frac{(\check{\cK}_\infty)_1(z)}{(\check{\cK}_0)_1(z)}=
    z\lim_{r\uparrow 1}\frac{-((I-r^2\cH^*\cH)^{-1}\cH^*\bar t)(z)}
    {((I-r^2\cH^*\cH)^{-1}1)(z)}
\end{equation}
 and $\psi_\cH$ is the outer function
\begin{equation}\label{psi}
  \psi_\cH(z)\psi_\cH(0)=\lim_{r\uparrow 1}\frac 1{((I-r^2\cH^*\cH)^{-1}1)(z)},\quad \psi_\cH(0)>0.
\end{equation}
 Moreover, $|\psi_\cH|^2+|\phi_\cH|^2=1$.
\end{theorem}

Similarly to \eqref{defphi}, by \eqref{rk8} we define
\begin{equation}\label{defphi2}
   \phi:= z\frac{({\cK}_\infty)_1(z)}{({\cK}_0)_1(z)}=\bar a_{-1}\frac{R(0)-R}{R(0)+R}
\end{equation}
and the outer function $\psi$, $\psi(0)>0$, such  that $|\psi|^2+|\phi|^2=1$.

\begin{remark}\label{rk3.4}
Let us note the important identity
\begin{equation}\label{verbk}
    \frac{({\cK}_\infty)_1(0)}{({\cK}_0)_1(0)}=-\frac{\bar a_{-1}} 2 R'(0)=
    -\bar a_{-1}\langle\fA^{-1}e_0,e_0\rangle=\langle \bar a_0e_0+\rho_0 e_1,e_0\rangle=\bar a_0.
\end{equation}
\end{remark}

As in \eqref{dcrp-nehclass} we  consider the collection of functions
\begin{equation}\label{aak10}
    f=-\frac{\psi}{\bar\psi}\frac{\cE+\bar\phi}{1+\phi\cE},\quad \cE\in H^\infty,\ \|\cE\|\le 1.
\end{equation}
Let us note $\psi/(1+\cE\phi)\in H^2$. Indeed,
$$
\left|\frac{\psi}{1+\cE\phi}\right|^2\le \frac{1-|\cE\phi|^2}{|1+\cE\phi|^2}=
\Re\frac{1-\cE\phi}{1+\cE\phi}\in L^1,
$$
and $1+\cE\phi$ is an outer function.
Therefore, due to the identity
$$
f=-\frac{\psi}{\bar\psi}\bar\phi
-\frac{\cE\psi^2}{1+\phi\cE},
$$
 all of them correspond to the same Hankel operator.

 \begin{definition}\cite{arov, AD1}
 A function $\phi$
 \begin{equation}\label{propphi}
    \phi\in H^\infty,\ \|\phi\|\le 1,\
    \phi(0)=0,\ \log(1-|\phi|^2)\in L^1,
 \end{equation}
 is called Arov-regular if the set \eqref{aak10} describes the collection of \textit{all} symbols with the same Hankel operator.

  It is called Arov-singular if $f_0=-\frac{\psi}{\bar\psi}\bar\phi\in H^\infty$. Equivalently, the entries of the unitary valued matrix function
 \begin{equation}\label{scatmat}
    \begin{bmatrix}f_0&\psi\\
    \psi&\phi
    \end{bmatrix}
 \end{equation}
 belong to $H^\infty$.
 \end{definition}
 It means that \eqref{aak10} with a singular $\phi$ describes a proper subclass of the Schur class, $\cH_f=0$ for all $f$. The Potapov-Ginzburg transform of the matrix \eqref{scatmat}
 \begin{equation}\label{pgt}
    A:=\frac 1{\psi}\begin{bmatrix}f_0\phi-\psi^2& f_0\\
    \phi& 1
    \end{bmatrix}=\begin{bmatrix}-\frac 1{\bar\psi}&-\frac{\bar\phi}{\bar\psi}\\
    \frac{\phi}{\psi}&\frac{1}{\psi}
    \end{bmatrix}
 \end{equation}
 is an (Arov-singular) $j$-inner matrix function,
 $$
 A(z)^*jA(z)-j\ge 0,\ z\in\bbD;\ A(t)^*j A(t)-j= 0,\ t\in\bbT, \ j=
  \begin{bmatrix}-1&0\\
    0& 1
    \end{bmatrix}.
 $$
 \begin{theorem}[Arov]  Every function of the form \eqref{propphi} possesses (Arov-) singular-regular factorization:
 \begin{equation}\label{asrf}
     \frac 1\psi\begin{bmatrix}\phi& 1
    \end{bmatrix}=\frac 1\psi_\cH\begin{bmatrix}\phi_\cH& 1
    \end{bmatrix}A.
 \end{equation}
 \end{theorem}

 \begin{theorem}\label{th3.6}\cite{arov}, see also \cite{Kh}.
 Let $\phi$ be a function of the form \eqref{propphi}, equivalently, let the Herglotz class function
 \begin{equation}\label{nado2}
    R(z)=\frac{1-a_{-1}\phi(z)}{1+a_{-1}\phi(z)}=\int\frac{t+z}{t-z}{d\sigma},\quad
    a_{-1}\in\bbT,
 \end{equation}
 be associated with the Szeg\"o class measure $\sigma$, $\sigma(\bbT)=1$.
 We set
 \begin{equation}\label{nado}
    \cH=\cH_{s},\ s=-a_{-1}\frac{D}{D_*},\ D(z)= \frac{\psi(z)}{1+a_{-1}\phi(z)}.
 \end{equation}
 The function $\phi$ is Arov-regular if and only if one of the following equivalent conditions hold:
 \begin{itemize}
\item[(i)]
$M(\sigma, a_{-1})=M_s(D)=\check M_s(\cH)$.
\item[(ii)] The reproducing kernel
 $\cK_0$ belongs to $\check M_s(\cH)$.
 \item[(iii)] $\phi=\phi_\cH,\quad \psi=\psi_{\cH}$.
 \item[(iv)] $\psi(0)=\psi_\cH(0)$, that is the following limit exists
\begin{equation}\label{lim8}
    \lim_{r\uparrow 1}\langle(I-r^2\cH_s^*\cH_s)^{-1}1,1\rangle=\frac{1}{D^2(0)}=\frac{1}{\psi^2(0)}.
\end{equation}
\end{itemize}
 \end{theorem}

\begin{definition}Let $\fA\in \Sz$ correspond to the spectral data $\sigma, a_{-1}$. We say that
$\fA$ is regular, $\fA\in \Sz^{reg}$, if the associated function
\begin{equation}\label{assphi}
    \phi(z)=\bar a_{-1}\frac{R(0)-R(z)}{R(0)+R(z)}, \quad R(z)=\int\frac{t+z}{t-z}d\sigma(t),
\end{equation}
is Arov-regular.
\end{definition}

\begin{theorem}
The correspondence $\fA\mapsto s$, $\fA\in \Sz_{a.c.}$, is one to one precisely on the subclass
 $\Sz^{reg}$.
\end{theorem}

 \begin{proof} Let $\phi$ be not regular. Then it is of the form \eqref{asrf}, where $A$ is a non constant $j$-inner matrix function  with the entries $\phi_A$, $\psi_A$.

 First we assume that the vector $A\begin{bmatrix} a_{-1}\\1\end{bmatrix}$ is collinear to a constant, that is
 \begin{equation}\label{const1}
    A(z)\begin{bmatrix} a_{-1}\\1\end{bmatrix}=\begin{bmatrix} c_1\\c_2\end{bmatrix}\omega(z),
    \ |c_1|^2+|c_2|^2\not=0.
 \end{equation}
We claim that in this case the associated $R$ function
 $$
 R(z)=\frac{1-a_{-1}\phi}{1+a_{-1}\phi}=\int\frac{t+z}{t-z}|D(t)|^2 dm+\int\frac{t+z}{t-z}d\sigma_s(t)
 $$
 has the absolutely continuous component with the density proportional to the density of a canonical
 measure
 \begin{equation}\label{const2}
    D=D_{\tilde\cH}\psi_A(0),\quad
 \tilde R(z)=\frac{1-a_{-1}\phi_{\tilde\cH}}{1+a_{-1}\phi_{\tilde\cH}}=\int\frac{t+z}{t-z}|D_{\tilde\cH}(t)|^2 dm,
 \end{equation}
 where $\phi_{\tilde\cH}=\bar a_{-1}\tilde a_{-1}\phi_\cH$,
  $|\tilde a_{-1}|=1$.
 But also it has a non-trivial singular component $\sigma_s(\bbT)=1-|\psi_A(0)|^2$, that is $\fA\not\in \Sz_{a.c.}$

 Let us show that, in fact, $\omega(z)= const$ in \eqref{const1}. We note that we have here a so called $j$-neutral vector, i.e.,
 $$
 \begin{bmatrix} \bar c_1& \bar c_2\end{bmatrix}j\begin{bmatrix} c_1\\c_2\end{bmatrix}=
 -|c_1|^2+|c_2|^2=0.
 $$
 Indeed, it follows from the fact that $A$ is $j$-unitary on the boundary
 \begin{equation*}\label{const32}
 \begin{split}
    (-|c_1|^2+|c_2|^2)|\omega(t)|^2=
    \begin{bmatrix} \bar a_{-1} &1\end{bmatrix}A(t)^* jA(t)\begin{bmatrix} a_{-1}\\1\end{bmatrix}\\
=\begin{bmatrix} \bar a_{-1}&1\end{bmatrix}j\begin{bmatrix} a_{-1}\\1\end{bmatrix}
 =-|a_{-1}|^2+1=0.
 \end{split}
 \end{equation*}
 Now we use the Schwartz Lemma for $j$-expanding matrix functions
 \begin{equation}\label{const4}
     \begin{bmatrix} \frac{A(z)jA^*(z)-j}{1-|z|^2}&\frac{A(z)-A(0)}{z}\\
   *  & {A^*(0)jA(0)-j}
     \end {bmatrix}\ge 0.
 \end{equation}
 Since
 $$
   \begin{bmatrix}  a_{-1} \\1\end{bmatrix}^*(A^*(0)jA(0)-j)\begin{bmatrix}  a_{-1} \\1\end{bmatrix} =
   (-|c_1|^2+|c_2|^2)|\omega(0)|^2-(-|a_{-1}|^2+1)=0
 $$
 we get from \eqref{const4}
 \begin{equation*}\label{const5}
 \begin{split}
 \begin{bmatrix} \frac{A(z)jA^*(z)-j}{1-|z|^2}&\frac{A(z)-A(0)}{z}\begin{bmatrix}  a_{-1} \\1\end{bmatrix}\\
   *  & \begin{bmatrix}  a_{-1} \\1\end{bmatrix}^*(A^*(0)jA(0)-j)\begin{bmatrix}  a_{-1} \\1\end{bmatrix}
     \end {bmatrix}&\\=
      \begin{bmatrix} \frac{A(z)jA^*(z)-j}{1-|z|^2}&\begin{bmatrix} c_1\\c_2\end{bmatrix}\frac{\omega(z)-\omega(0)}{z}\\
   *& 0
     \end{bmatrix}\ge 0&,
     \end{split}
 \end{equation*}
 which implies $\omega(z)-\omega(0)=0$.

 Therefore we proved that
 \begin{equation}\label{const5}
    A(z)\begin{bmatrix} a_{-1}\\1\end{bmatrix}=A(0)\begin{bmatrix} a_{-1}\\1\end{bmatrix}=
    \begin{bmatrix} \tilde a_{-1}\\1\end{bmatrix}\frac 1{\psi_A(0)}
 \end{equation}
 with a certain $\tilde a_{-1}\in \bbT$. By \eqref{asrf} we get
 \begin{equation}\label{const6}
 \begin{split}
 \frac 1 D=&
    \frac 1\psi \begin{bmatrix} \phi &1\end{bmatrix}\begin{bmatrix} a_{-1}\\1\end{bmatrix}=
     \frac 1\psi_\cH \begin{bmatrix} \phi_\cH &1\end{bmatrix}A\begin{bmatrix} a_{-1}\\1\end{bmatrix}\\
    =&\frac 1\psi_\cH \begin{bmatrix} \phi_\cH &1\end{bmatrix}\begin{bmatrix} \tilde a_{-1}\\1\end{bmatrix}
    \frac 1 {\psi_A(0)}=\frac 1{\psi_A(0)D_{\tilde\cH}}.
    \end{split}
 \end{equation}
Thus \eqref{const2} is proved.

 It remains to consider the case
 \begin{equation}\label{nonconst1}
    A(z)\begin{bmatrix} a_{-1}\\1\end{bmatrix}=\begin{bmatrix} \cE(z)\\ 1\end{bmatrix}\omega(z),
 \end{equation}
 where the inner function $\cE$ is not a constant, that is,
 the symbol $s$   has a nontrivial inner function $\cE$  in its canonical representation \eqref{dcrp-nehclass}. We have to show  that
 there are at least two different CMV matrices from $\Sz_{a.c.}$ with the given scattering function.

For $t_1, t_2\in \bbT$, $t_1\not= t_2$, let $e^{ic_k}$ be such that
 $e^{ic_k}(1+\bar t_k\cE(0))>0$. We define
 $$
 D_k=C_k e^{ic_k}(1+\bar t_k\cE)\frac{\psi_\cH}{1+\cE\phi_{\cH}},\quad
 (a_{-1})_k=e^{-2ic_k}t_k,
 $$
 where the normalizing positive constants $C_k$ are chosen from the conditions
 $\int|D_k|^2 dm=1$.
 Both data $\{|D_k|^2 dm,(a_{-1})_k\}$, $k=1,2$, produce the same function $s$ corresponding to the given $\cE$.
 \end{proof}

Here is a sufficient condition known in the context of the AAK theory.
\begin{theorem}
Let $\fA\in \Sz$ and let  its spectral measure be absolutely continuous, $d\sigma =w dm$. If $1/w\in L^1$ then $\fA\in \Sz^{reg}$.
\end{theorem}
\begin{proof} For the given scattering function $s$ we have the canonical representation
$$
s=-a_{-1}\frac{D}{\bar D}=-\cE\frac{\psi_\cH}{\bar\psi_\cH}\frac{1+\bar\phi_\cH  \bar \cE}{1+\phi_\cH\cE}.
$$
Therefore,
\begin{equation}\label{h1h1}
    G:=\frac 1{D}\cE\frac{\psi_\cH}{1+\phi_\cH\cE}=a_{-1}\overline{\frac 1{D}\frac{\psi_\cH}{1+\phi_\cH\cE}}.
\end{equation}
The function $\frac{\psi_\cH}{1+\phi_\cH\cE}$ belongs to $H^2$ and, due to the assumption, $\frac 1 D\in H^2$. Thus $G$ belongs to $H^1$, that is, all its negative Fourier coefficients vanish. From the second representation $\bar G\in H^1$. That is, all positive  Fourier coefficients of $G$ vanish. Therefore $G$ is constant. Since $\cE$ is the inner part of $G$, we have $\cE=const$. Using the normalization
$D(0)>0, \psi_\cH(0)>0$, we get $\cE=a_{-1}$.
\end{proof}

A similar sufficient condition $1/\psi\in H^2$ is given in Proposition \ref{prop4.2}. For a weaker condition on $|\psi|$, which ensures regularity of $\phi$, see \cite{VYu2}.

\section{Helson-Szeg\"o Theorem and Boundedness of the GLM Transform}

We define the orthonormal system in $M(\sigma,a_{-1})$
$$
\fe_n(z):=\cF(z,P_n).
$$
Recall, the first component $(\fe_n)_1(z)$ is holomorphic for $|z|<1$ and the second,
$(\fe_n)_2(z)$, for $|z|>1$. Moreover,
due to the orthogonality property of the Laurent  polynomials $P_n$ we have
\begin{equation}\label{hs1}
\begin{split}
\fe_{2n}=&\left(\frac{1}{D(0)\pi^{(2n)}_{2n}}z^{n}+\dots, O(\frac 1{z^{n+1}})\right),\\
\fe_{2n+1}=&\left(O(z^{n+1}),\frac{1}{D(0)\pi^{(2n+1)}_{2n+1}}\frac 1{z^{n+1}}+\dots \right).
    \end{split}
\end{equation}

\begin{definition}
The (lower--triangular) matrix of the  Gelfand-Levitan-Marchenko (GLM) transformation is defined by the following relation:

\begin{equation}\label{hs2}
\begin{bmatrix}\fe_0&\fe_1&\fe_2&\dots \end{bmatrix}=
\begin{bmatrix}\begin{bmatrix}1\\0\end{bmatrix}&\begin{bmatrix}0\\1/z\end{bmatrix}&
\begin{bmatrix}z\\0\end{bmatrix}&\dots\end{bmatrix}\cM,
\end{equation}
that is,
\begin{equation}\label{defms}
\begin{split}
   \fe_{2n}=& \cM^{2n}_{2n}\begin{bmatrix}z^n\\ 0\end{bmatrix}+\cM^{2n}_{2n+1}\begin{bmatrix}0\\ 1/z^{n+1}\end{bmatrix}+
   \cM^{2n}_{2n+2}\begin{bmatrix}z^{n+1}\\ 0\end{bmatrix}+\dots\\
     \fe_{2n+1}=& \cM^{2n+1}_{2n+1}\begin{bmatrix}0\\ 1/z^{n+1}\end{bmatrix}+\cM^{2n+1}_{2n+2}\begin{bmatrix}z^{n+1}\\ 0\end{bmatrix}+
   \cM^{2n+1}_{2n+2}\begin{bmatrix} 0\\ 1/z^{n+2}\end{bmatrix}+\dots
   \end{split}
\end{equation}
\end{definition}
Let us point out that generally only matrix elements of $\cM$ are well defined.
\begin{proposition}\label{prop4.2}
The vector $\cM e_0$ belongs to $l^2(\bbZ_+)$ if and only if $1/\psi\in H^2$. In particular,
$\cM e_0\in l^2(\bbZ_+)$ implies $\fA\in \Sz^{reg}$.
\end{proposition}

\begin{proof} By the definition, see \eqref{rk8},
$$
\fe_0=\left(\frac 1{\psi\psi(0)},\frac{\bar\phi}{\bar\psi\psi(0)}\right).
$$
Therefore $\sum_{j\ge 0}|\cM^0_{2j}|^2<\infty$ is equivalent to $\frac 1{\psi\psi(0)}\in H^2$ and
$\sum_{j\ge 0}|\cM^0_{2j+1}|^2<\infty$ is equivalent to $\frac {\bar\phi}{\bar\psi\psi(0)}\in H^2_-$.

Let us chose $f_0=-\bar\phi\frac{\psi}{\bar\psi}$ as a symbol of the associated Hankel operator $\cH$.
Since $1/\psi\in H^2$, we get
$$
(I-\cH^*\cH)\frac{1}{\psi\psi(0)}= \frac{1}{\psi\psi(0)}-P_+(\phi\frac{\bar\psi}{\psi})\frac{\bar\phi}{\bar \psi\psi(0)}=\frac{1}{\psi\psi(0)}-P_+(1-|\psi|^2)\frac{1}{ \psi\psi(0)}=1.
$$
Therefore, $\psi=\psi_\cH$. By Theorem \ref{th3.6} $\fA\in \Sz^{reg}$.
\end{proof}

\begin{theorem}\label{th4.3} Let $\fA\in \Sz^{reg}$, that is, $M(\sigma, a_{-1})=\check M_s$. In this case the orthonormal basis \eqref{hs1} is of the form
\begin{equation}\label{onmb}
\fe_{2n}=\begin{bmatrix}t^n&0\\ 0& \bar t^n\end{bmatrix}\frac{\check \cK^{st^{2n}}_{0}}
{\|\check \cK^{st^{2n}}_{0}\|},
\quad
\fe_{2n+1}=- a_{-1}
\begin{bmatrix}t^{n+1}&0\\ 0 & \bar t^n\end{bmatrix}\frac{\check \cK^{st^{2n+1}}_{\infty}}
{\|\check \cK^{st^{2n+1}}_{\infty}\|},
\end{equation}
where $\check \cK^{st^{n}}_{0}$, $\check \cK^{st^{n}}_{\infty}$ are the reproducing kernels \eqref{defchk} of $\check M_{st^{n}}$. In particular,
\begin{equation}\label{rhokern}
\begin{split}
    \cM_{2n}^{2n}=&\sqrt{(\check \cK^{st^{2n}}_{0})_1(0)}=
    \sqrt{\lim_{r\uparrow 1}
    \langle(I-r^2\cH_{st^{2n}}^*\cH_{st^{2n}})^{-1}1,1\rangle}\\
    =&
    \frac{\rho_0\rho_1\dots\rho_{2n-1}}{D(0)}=\frac 1{\rho_{2n}\rho_{2n+1}\dots},\\
    \cM_{2n+1}^{2n+1}=&-a_{-1}\sqrt{(t\check \cK^{st^{2n+1}}_{\infty})_2(\infty)}=
    -a_{-1}\sqrt{\lim_{r\uparrow 1}
    \langle(I-r^2\cH_{st^{2n+1}}\cH_{st^{2n+1}}^*)^{-1}\bar t,\bar t\rangle}\\=&
    -a_{-1}\frac{\rho_0\rho_1\dots\rho_{2n}}{D(0)}=\frac {-a_{-1}}{\rho_{2n+1}\rho_{2n+2}\dots}.
    \end{split}
\end{equation}

\end{theorem}

Recall that the space $\check M_{st^{n}}$, in fact depends on the Hankel operator $\cH_{st^{n}}$,
that is, of the negative Fourier coefficients of the function ${st^{n}}$.

\begin{remark} Theorem \ref{th4.3} provides an  algorithm for inverse scattering, indeed, \eqref{onmb} gives also an explicit formula for the Verblunsky coefficients, see Remark \ref{rk3.4}:
\begin{equation}\label{expla}
    \bar a_n= \frac{(\check\cK^{st^n}_\infty)_1(0)}{(\check\cK^{st^n}_0)_1(0)}.
\end{equation}
\end{remark}

\begin{proof}[Proof of Theorem \ref{th4.3}] Let us note that $(\fe_n)_1(0)=0$ for $n\ge 1$. That is all vectors spanned by this system are orthogonal to $\check \cK^{s}_{0}$. Therefore $\check \cK^{s}_{0}$ is collinear to $\fe_0$ and coincides with the initial basis vector after the appropriate normalization. A dense set in the orthogonal complement in $\check M_{s}$ is of the form $\cF=(tF_1,F_2)$, $F_1\in H^2$, $F_2\in H^2_-$. Now,
by  definition $(tF_1,F_2)\in\check M_{s}$ means that $(F_1,F_2)\in\check M_{ts}$. Thus we have the decomposition
$$
\check M_{s}=\{\check \cK^{s}_{0}\}\oplus\begin{bmatrix}t&0\\ 0 & 1\end{bmatrix}\check M_{ts}.
$$
Continue in this way we get \eqref{onmb}.

To prove \eqref{rhokern} we use \eqref{hs1} and \eqref{pirho}.
\end{proof}

\begin{theorem}
Let $\fA\in \Sz^{reg}$. Then the following upper-lower triangular factorization holds true
\begin{equation}\label{upl1}
    \lim_{r\uparrow 1}U^*\begin{bmatrix} I&r\cH^*\\
    r\cH& I
    \end{bmatrix}^{-1}U=\cM \cM^*,
\end{equation}
where $U:l^2(\bbZ_+)\to H^2\oplus H^2_-$ is the reordering of the standard basis
$$
U e_{2n}=t^n\oplus 0,\quad U e_{2n+1}=0\oplus  t^{-n-1}.
$$
In particular, $\|\cM\|<\infty$ if and only if $\|\cH\|<1$.
\end{theorem}

\begin{proof} Let $\cM_r$ be the GLM transform that corresponds to the Hankel operator $r\cH$.
In this case directly from \eqref{hs2} we have
$$
\cM_r^*U^*\begin{bmatrix} I&r\cH^*\\
    r\cH& I
    \end{bmatrix}U\cM_r =I.
    $$
    We rewrite this into the form
    \begin{equation*}\label{upl2}
    U^*\begin{bmatrix} I&r\cH^*\\
    r\cH& I
    \end{bmatrix}^{-1}U=\cM_r \cM_r^*.
\end{equation*}
We can pass here to the limit since only a finite number of elements of $\cM_r$ are involved in the entry of the product (recall it is a lower triangular matrix) and $(\cM_r)^k_l\to \cM^k_l$ for fixed $k$ and $l$.
\end{proof}

\begin{definition}
We say that $\fA$ belongs to the Helson-Szeg\"o class, $\fA\in \HSz$ if the matrix of the GLM transform generates a bounded operator in $l^2(\bbZ_+)$, $\|\cM\|<\infty$.
\end{definition}

\begin{proposition}\label{prop4.5}
The GLM transformation is a bounded operator in $l^2(\bbZ_+)$ if and only if $\sigma$ is absolutely continuous and the Riesz projections $P_\pm$ are bounded operators in $L^2_{w^{-1}dm}$.
\end{proposition}

\begin{proof} It was proved in Proposition \ref{prop4.2} that $\sigma$ is absolutely continuous and $M_s(D)=\check M_s$. Let
$$
\cF=(F_1,F_2)=\sum  \tilde f_{k}\fe_k, \ \ \tilde f\in l^2(\bbZ_+).
$$
Then
$$
\langle\cM\tilde f,\cM\tilde f\rangle\le C\langle\tilde f,\tilde f\rangle
$$
means
\begin{equation}\label{hsineq}
    \|F_1\|^2+\|F_2\|^2\le C\|\cF\|^2_s,
\end{equation}
where in the LHS we have the standard $L^2$ norm.

On the other hand, according to definition \eqref{toscat1}
$$
F_1=\frac 1 D P_+ g,\ \ F_2=-\frac {a_{-1}}{D_*} P_- g, \ g=wf.
$$
We put these in \eqref{hsineq}. Due to $\|\cF\|^2_s=\|f\|^2_{L^2_{wdm}}$, we get
$$
\langle w^{-1}P_+ g, P_+g\rangle+\langle w^{-1}P_-g, P_-g\rangle\le
C\langle w f,f\rangle=C\langle w^{-1} g,g\rangle.
$$
\end{proof}

Recall that a weight $w$ satisfies $A_2$ (or Hunt-Muckenhoupt-Wheeden) condition
if for all arcs $I\subset \bbT$ the following supremum is finite
\begin{equation}\label{hmw}
    \sup_{I}\langle w\rangle_I\langle w^{-1}\rangle_I<\infty,
\end{equation}
where  $\langle w\rangle_I=\frac 1{|I|}\int_{I}wdm$.

We combine the classical   Helson-Szeg\"o and Hunt-Muchenhoupt-Wheeden Theorems, see e.g. \cite{Nik},
with Proposition \ref{prop4.5}.

\begin{theorem}\label{th4.7} The following are equivalent:
\begin{itemize}
\item[(i)] $\fA\in \HSz$,
\item[(ii)] $w\in A_2$,
\item[(iii)] $\|\cH_s\|<1$, $w^{-1}\in L^1$,
\item[(iv)] $s\in HS$.
\end{itemize}
\end{theorem}

\section{Golinskii-Ibragimov Theorem and Faddeev-Marchenko type scattering theorem}

The following theorem was proved in \cite{IR}, see also \cite{pel} and \cite{KhP}, where the general case $B_p^{1/p}$ was considered.

\begin{theorem}
Let $s(t)$ be an unimodular function of the class $B_2^{1/2}$. Then there exists a unique representation
\begin{equation}\label{pelth}
    s(t)=t^N e^{iv(t)}
\end{equation}
where $v(t)\in B_2^{1/2}$.
\end{theorem}
The integer $N$ in the representation \eqref{pelth} is called the index of $s(t)$. It can be computed as the winding number of the harmonic extension $s(rt)$ in the unit disk  of the  given function   for $r$ sufficiently close to $1$.

\begin{definition} We say that $\fA$ belongs to the Golinskii-Ibragimov class, $\fA\in \GI$ if
\begin{equation}\label{defgi}
    \sum k|a_k|^2<\infty
\end{equation}
\end{definition}

\begin{theorem}\label{th5.3}
A unimodular function $s(t)$ is the scattering function for a unique $\fA\in \GI$ if and only if $s(t)\in B_2^{1/2}$ and its index $N=0$. Moreover
$\GI\subset \HSz$.
\end{theorem}

 Note, statements like NLFT belongs to a certain class if and only if LFT belongs to the same class  is typical for the Faddeev--Marchenko  theory, see e.g.  \cite[Theorem 3.3.3]{Mar1}, or the newest results of this type \cite{EMT2}.

Theorem \ref{th5.3}  is the scattering counterpart  of the (spectral) Golinskii-Ibragimov version of the Strong Szeg\"o Limit Theorem (see \cite{GI} and \cite{Sv1}).

\begin{theorem}\label{th5.4} A measure $d\sigma =w dm+d\sigma_{s}$ is the spectral measure of $\fA\in \GI$ if and only if $\log w(t)\in B_2^{1/2}$ and  $\sigma_s=0$.
\end{theorem}

\begin{proof}[Proof of Theorem \ref{th5.3}]
Let $\fA\in \GI$. By Theorem \ref{th5.4}
$$
D=e^{\frac{u+i\tilde u}{2}},\ u:=\log w.
$$
Therefore
$$
s=-a_{-1} e^{i\tilde u}\in B_2^{1/2}.
$$

Conversely, let $s=e^{iv}\in B_2^{1/2}$. We define
$$
a_{-1}=-e^{i\int v dm}, \ w=Ce^{-\tilde v},\ C>0: \int w dm=1.
$$
Therefore $\fA$ with the spectral data $(wdm,a_{-1})$ belongs to $\GI$.

To show that this is the only CMV matrix of the Szeg\"o class $\Sz_{a.c.}$, which corresponds to the given scattering function we note that $P_-s\in B_2^{1/2}$ means precisely that
$\cH_s^*\cH_s$ belongs to the trace class. Therefor we have the alternative: 1)$\|\cH_s\|<1$, or, 2) there exists $g\in H^2$, $g\not=0$, such  that
$$
(I-\cH_s^*\cH_s)g=0.
$$
In the second case we have
$$
\|P_+sg\|=\|sg\|^2-\|P_-sg \|^2=0.
$$
Thus $h:=sg\in H_-^2$, and we get
$$
-a_{-1}Dg=\bar D h.
$$
The LHS is in $H^1$ and the RHS has all nonnegative Fourier coefficients equal to zero. Therefor
$Dg=0$, but this contradicts to  $g\not=0$.

 Thus $\|\cH_s\|<1$ and $w^{-1}=1/C e^{\tilde v}\in L^1$. By Theorem \ref{th4.7} $\fA\in \HSz$ which guarantees uniqueness of the inversion problem.
\end{proof}

\begin{remark}\label{rk5.5}
To our best knowledge scattering for CMV matrices was not studied in the frame of the standard Faddeev-Marchenko approach, that is we cannot compare Theorem \ref{th5.3} with a "traditional" one.
 But we can discuss certain conditions related to coefficients of Jacobi matrices.   The first one is the classical  scattering theory condition, see \cite{GC, Gu} and also \cite{T},
\begin{equation}\label{guscond}
    \sum n(|p_n-1|+|q_n|)<\infty.
\end{equation}
The second condition was obtain by E. Ryckman \cite{Ry} who proved a counterpart of the Strong Szeg\"o Theorem for Jacobi matrices.  It corresponds to   $\sum n|a_n|^2<\infty$ in the CMV case and
 is of the form
 \begin{equation}\label{rycond}
    \sum_{n=1}^\infty n\{|\sum_{k=n}^\infty(p_k-1)|^2+|\sum_{k=n}^\infty q_k|^2\}<\infty.
\end{equation}
Evidently, \eqref{guscond} implies \eqref{rycond} (we discuss only the behavior of the $q_k$'s) :
\begin{equation*}
    \begin{split}
     \sum_{n=1}^\infty n|\sum_{k=n}^\infty q_k|^2\le& \sum_{n=1}^\infty \sum_{k=n}^\infty(\sum_{l=n}^\infty n|q_l|)|q_k|
     =\sum_{k=1}^\infty \sum_{n=1}^k(\sum_{l=n}^\infty n|q_l|)|q_k|\\
     \le&\sum_{k=1}^\infty \sum_{n=1}^k(\sum_{l=1}^\infty l|q_l|)|q_k|=
     (\sum_{l=1}^\infty l|q_l|)(\sum_{k=1}^\infty k|q_k|).
    \end{split}
\end{equation*}
Note that for a typical example $q_n=\frac 1{n^\beta}$ both conditions \eqref{guscond},  \eqref{rycond}
require $\beta>2$. However, for the oscillating sequence  $q_n=\frac {(-1)^n}{n^\beta}$ we have
$\beta>2$ for \eqref{guscond} and $\beta>1$ for \eqref{rycond}.
\end{remark}
\begin{remark} Note that for a typical $A_2$ weight $w(t)=|t-1|^{2\gamma_1}|t+1|^{2\gamma_2}$,
$\gamma_k>-1/2$,
the so called Jacobi (OPUC) Verblunsky coefficients are of the form
\begin{equation*}\label{golinskii}
    a_n=-\frac{\gamma_1-(-1)^n\gamma_2}{n+1+\gamma_1+\gamma_2},
\end{equation*}
which violates \eqref{defgi}. That is, a characterization of the class $\HSz$ in terms of the Verblunsky coefficients looks like an extremely interesting challenging problem.
\end{remark}

\section{Widom's Formula}
The so-called Widom Formula has a natural proof in the frame of the scattering theory. Roughly speaking it says, see \eqref{upl1} and \eqref{rhokern},
\begin{equation}\label{rwidf}
    \det\begin{bmatrix}I&\cH^*\\ \cH& I
    \end{bmatrix}^{-1}=\prod_{n=0}^\infty |\cM^n_n|^2=\frac 1{\rho_0^2\rho_1^4\rho_2^6\dots}.
\end{equation}
Here the determinant has sense if $\cH$ is of the trace class. But, in fact, a stronger statement holds true.
\begin{theorem}\label{th6.1} Let $\cH$ be of the Hilbert-Schmidt class. Then
\begin{equation}\label{widf}
    \det(I-\cH^*\cH)=\rho_0^2\rho_1^4\rho_2^6\dots
\end{equation}
\end{theorem}

To get \eqref{widf} we prove a counterpart of the factorization formula
\eqref{upl1} for the matrix $\lim_{r\uparrow1}(I-r^2\cH^*\cH)^{-1}$. Consider a subspace of
$\check M_s$ consisting of the vectors of the form
$$
\text{clos}\{\cF=(F,-\cH F): F\in H^2\}\subset \check M_s.
$$
In fact, $\cF$ belongs to this subspace if and only if $F$ belongs to
$$
\check M_s^+=\text{clos}\{F\in H^2:\|F\|_s^2=\langle(I-\cH^*\cH) F, F\rangle\}.
$$
By the way, we have the
orthogonal decomposition
$$
\check M_s=\check M_s^+\oplus H^2_-
$$
in the following  sense
$$
\cF=(F,-\cH F)\oplus (0,G), \quad F\in \check M_s^+,\ G\in H^2_-.
$$
Also,
$$
\check \cK^{s}_{0}=((\check \cK^{s}_{0})_1,-\cH(\check \cK^{s}_{0})_1),
$$
 and therefore $(\check \cK^{s}_{0})_1\in \check M_s^+$ is the reproducing kernel in this subspace.

Similar to Theorem \ref{th4.3} we have
\begin{theorem} Let $\fA\in \Sz^{reg}$. The system of vectors
\begin{equation}\label{onmbplus}
\ff_{n}=\frac{t^n(\check \cK^{st^{n}}_{0})_1}
{\|\check \cK^{st^{n}}_{0}\|}
\end{equation}
forms an orthonormal basis in $\check M_s^+$.
\end{theorem}

Similar to \eqref{hs2} we define
the (lower--triangular) matrix
\begin{equation}\label{wfl}
\begin{bmatrix}\ff_0&\ff_1&\ff_2&\dots \end{bmatrix}=
\begin{bmatrix}1& z & z^2\dots\end{bmatrix}\cL.
\end{equation}
In this case, similar to \eqref{upl1}, we have
\begin{equation}\label{wup}
    \lim_{r\uparrow1}(I-r^2\cH^*\cH)^{-1}=\cL\cL^*.
\end{equation}

Finally, we note that
$$
 \langle(I-r^2\cH_{st^{2n+1}}\cH_{st^{2n+1}}^*)^{-1}\bar t,\bar t\rangle
=
\langle(I-r^2\cH_{st^{2n+1}}^*\cH_{st^{2n+1}})^{-1}1,1\rangle
$$
and therefore, by \eqref{rhokern}, \eqref{onmbplus} and \eqref{wfl}, we have
\begin{equation}\label{lnn}
    \cL^n_n=\sqrt{\langle(I-r^2\cH_{st^{n}}^*\cH_{st^{n}})^{-1}1,1\rangle}=
\frac 1{\rho_{n}\rho_{n+1}\dots}
\end{equation}
for both even and odd $n$'s.

 \begin{proof}[Proof of Theorem \ref{th6.1}] Since $\|\cH\|<1$ and $\tr(\cH^*\cH)<\infty$ we have
 $$
 (I-\cH^*\cH)^{-1}=I+\Delta,
 $$
 where $\Delta\ge 0$, $\tr\Delta<\infty$. Let $\Delta^{(n)}$ be the initial $n\times n$ block of the matrix $\Delta$. Then
 $$
 \det(I+\Delta)=\lim_{n\to\infty}\det(I+\Delta^{(n)}).
 $$
 Due to  the triangular factorization \eqref{wup}
 $$
 \det(I+\Delta^{(n)})=\prod_{k=0}^{n-1} |\cL^k_k|^2,
 $$
  where $\cL^k_k$ is given by \eqref{lnn}.
 \end{proof}

\bibliographystyle{amsplain}

\bigskip

{\em Mathematics Division, Institute for Low Temperature Physics and

Engineering, 47 Lenin ave., Kharkov 61103, Ukraine}

{E-mail: leonid.golinskii@gmail.com}

\bigskip
{\em Department of Mathematics, University of Massachusetts Lowell,

Lowell, MA, 01854, USA}

E-mail: Alexander\_Kheifets@uml.edu

\bigskip
{\em Institute for Analysis, Johannes Kepler University of Linz,

A-4040 Linz, Austria}

E-mail: Franz.Peherstorfer@jku.at

E-mail: Petro.Yudytskiy@jku.at

\end{document}